\documentclass{amsart}

\newcounter{secti}

\newcounter{item}[section]
\renewcommand{\theitem}{\thesection.\arabic{item}}
\newcommand{\qitem}[2]{\refstepcounter{item}
  {\bf \theitem\  #1.}\  {#2 }}
\newcommand{\eitem}{\medskip}

\newcommand{\de}{\partial}
\newcommand{\der}[2]{\frac{\ddi #1}{\ddi #2}}
\newcommand{\pder}[2]{\frac{\de #1}{\de #2}}
\newcommand{\tends}{\rightarrow}
\newcommand{\di}{{\: \rm d}}   
\newcommand{\ddi}{{\rm d}}     
\newcommand{\E}{{\rm E\, }} 
\newcommand{\expd}{{\rm Exp\, }} 
\newcommand{\eop}{$ \Box $} 
\newcommand{\sign}{{\rm sgn\: } }

\newcommand{\be}{\begin{equation} }
\newcommand{\ee}{\end{equation}}
\newcommand{\beq}{\begin{eqnarray} }
\newcommand{\eeq}{\end{eqnarray}}

\newcommand{\eps}{\varepsilon}
\newcommand{\fai}{\raisebox{.1em}{$\varphi$}}

\newcommand{\cL}{{\mathcal L}}
\newcommand{\bx}{{\bf x}}
\newcommand{\by}{{\bf y}}
\newcommand{\bm}{{\bf m}}
\newcommand{\bn}{{\bf n}}
\newcommand{\R}{{\bf R}}
\newcommand{\cmn}{{\it CM}}
\newcommand{\cm}{{\it CM} }
\newcommand{\cmf}{\cm function }
\newcommand{\cmt}{\cmn. }
\newcommand{\cmz}{\cmn, }
\newcommand{\id}{{\it ID}}

\begin{document}

\begin{center}

{\Large\bf Completely monotone functions - a digest}

\medskip

{\large \em Milan Merkle}

\bigskip

\parbox{25cc}{
{\bf Abstract. }{\small This work has a purpose  to collect selected facts about the completely monotone ({\it CM}\/) functions that can
be found in books and papers devoted to
different areas of mathematics. We opted for lesser known ones, and for those which may help
determining whether or not a given function is completely monotone. In particular, we emphasize  the role of representation
of a {\it CM} function as the Laplace transform of a measure, and we present and discuss a little known
connection with log-convexity. Some of presented methods are illustrated by
several examples involving Gamma and related functions.

}
}

\end{center}

\medskip




\section{Introduction}
\label{intro}

A positive function defined on $(0,+\infty)$ of the class $C^{\infty}$, such that the sequence of its derivatives alternates signs
at every point, is called \index{completely monotone function} completely monotone (\cmn). A brief search in MathSciNet reveals total of 286 items that mention
this class of functions in the title from 1932 till the end of the year 2011; 98 of them have been published since the beginning of 2006.

This vintage topic was developed in 1920's/30's by
S. Bernstein, F. Hausdorff and V. Widder, originally with relation to so called moment problem, cf. \cite{bern29,haus21,haus23,wid31,wid34}.
 The much cited (but perhaps not that much read) Widder's book \cite{wid} contains a detailed account on properties
of \cm functions and their characterizations. The second volume of Feller's probability book \cite{fell2} discusses \cm functions through
their relationship with infinitely divisible measures, which are fundamental in defining  L\'evy processes. In past several decades,
L\'evy processes have gained popularity in financial models, as well as in biology and physics; this is probably a reason for increased interest
in \cm functions, too. There are also other interesting topics in Probability and Statistics where \cm functions play a role, see \cite{kimber74}
for one such topic. Aside from probability and measure theory, \cm and related functions appear in the field of approximations of functions, as
documented in the book \cite{meshf07} of 2007. Finally, they are naturally linked to various inequalities; several general
inequalities for \cm functions
can be found in \cite{marolk}, for a quite recent contribution in this area see \cite{auden2012}.

This text has a purpose  to collect well known  facts about the \cm functions, together with some less known ones, which may help
determining whether or not a given function is completely monotone. In that sense, this work can be thought of as being an extension and
supplement to another paper in the same spirit -- \cite{milsam01} by Miller and Samko. In particular, we emphasize  the role of representation
of a \cmf as a Laplace transform of a measure, and we present and discuss a little known (and even  less being used)
connection between \cmf and log-convexity. Some of methods discussed in sections \ref{reprcm}--\ref{inver} are illustrated by
several examples involving Gamma and related functions in  Section \ref{examp}. References and examples  reflect author's
preferences, and are by no means complete; the same can be said for the selection of topics that are discussed in this work.

\section{Representations of completely monotone functions}
\label{reprcm}

We start with a classical definition of \cm functions, and we present two  possible representations in terms of integral
transforms of measures and alternative representations for Steltjes transforms and \cm probability densities.

\qitem{Definition}{\label{defcm} A function $f$ defined on $(0,+\infty)$ is completely monotone if it has derivatives
of all orders and
\be
\label{cmdef}
 (-1)^kf^{(k)}(t) >  0,\qquad t\in (0,+\infty), k=0,1,2,\ldots \qquad\Box
 \ee
 }
In particular, this implies that each \cm
 on $(0,+\infty)$ is positive, decreasing and convex, with concave first derivative.

 \eitem

\qitem{Limit properties}{\label{limpro} By (\ref{cmdef}), there exist limits of $f^{(k)}(x)$ as $x\tends 0$ for any $k\geq 0$;
if those limits are finite, then $f$ can be extended
to $[0,+\infty)$ and (\ref{cmdef}) will also hold for $x=0$ (with strict inequality for all $k$). Limits at zero need not be finite, as
in $f(x)=\frac{1}{x}$, for example.

Clearly, $\lim_{x\tends +\infty} f^{(k)}(x)=0$ for all $k\geq 1$. The limit of $f(x)$ at $+\infty$ must be finite, and if it is
non-zero, then it has to be positive (for example, $f(x)=1+e^{-x}$).
}

\eitem

\qitem{Lemma}{\em \label{laptran} The function $f$ is \cm  if and only if
\cite{wid}
\be
\label{laptra}
 f(x) = \int_{[0,+\infty)} e^{-xt} \di \mu (t) ,
 \ee
where $\mu(t)$ is a positive measure on Borel sets of $[0,+\infty)$ (that is, $\mu (B)\geq 0$ for every Borel set $B\in \R_{+}$)
and the integral converges for
$0<x<+\infty$.}

 In other words, completely monotone functions are real one-side Laplace transforms
of a positive measure on $[0,+\infty)$. If the measure $\mu$ has an atom at $t=0$, then $\lim_{x\tends +\infty} f(x) >0$.
The measure $\mu$ is a probability measure if and only if
$\lim_{x\tends 0_{+}}f(0)=1$ (by monotone convergence theorem).

\eitem

The Lebesgue integral in (\ref{laptra})  can be expressed as
a Lebesgue-Stieltjes integral
\be
\label{laptras}
 f(x) = \int_{[0,+\infty)} e^{-xt} \di g (t) ,
 \ee
where $g(t)=\mu ([0,t])$ is the distribution function of $\mu$, with $g(0_{-})=0$.   For a
positive measure $\mu$, the function $g$ is non-decreasing, and by change of variables $t=-\log s$ we get

\qitem{Lemma}{\em \label{atran} The function $f$ is completely monotone on $(0,+\infty)$ if and only if
\be
\label{laptras2}
f(x) = \int_{[0,1]} s^x \di h(s),
\ee
where $h(s)=-g(-\log s)$ is a non-decreasing function.}

\eitem

If $f$ is a \cm which is the  Laplace transform of a measure $\mu$, as in (\ref{laptra}),
we write $f=\cL (\di\mu)$ or $f(x)=\cL (\di \mu(t))$. Similarly, the relation (\ref{laptras}) between $f$ and a distribution function $g$,
can be denoted as $f=\cL (\di g)$. If $\mu$ has a density $h$ with respect to Lebesgue measure, we write $f(x)=\cL (h(t)\di t)$ or only
$f=\cL (h)$. It follows from inversion formulae that each \cm  $f$ determines one positive measure $\mu$ via relation $f=\cL (\di \mu)$
and it is of  interest in many applications to find that measure.

\qitem{Remark}{\label{rmden} Since measures are determined by their Laplace transforms, if $f=\cL (\di \mu)$, then $f$ is \cm if and only if $\mu $ is
a positive measure. If there exists a continuous density $h$ of $\mu$, then $f$ is \cm if and only if $h(t)\geq 0$ for all $t\geq 0$.\eop
}

\eitem

Let us now observe a subclass of \cm functions which contains all functions $f$ that can be represented as  Stieltjes transform of some
positive measure $\mu $, that is,
\be
\label{stitra}
f(x) = \int_{[0,+\infty)}  \frac{\di \mu (s)}{x+s}
\ee

It is easy to verify that each function of the form (\ref{stitra}) with a positive measure $\mu$ is  \cmz hence $f=\cL (\nu)$, where
$\nu$ is a positive measure. To find $\nu$,  we start with
\[ \frac{1}{x+s} = \int_{[0,+\infty)} e^{-(x+s)u} \di u ,\]
and, after a change of order of integration we arrive at the following  result.

\qitem{Lemma}{{\em \label{stiltra} The Stieltjes transform of a positive measure $\mu$ as defined by (\ref{stitra})
can be represented as a Laplace
transform}
\[ f(x) = \int_{[0,+\infty)}e^{-xu} \left(\int_{[0,+\infty)} e^{-su}\di \mu(s)\right) \di u .\]

That is, $f=\cL (\nu)$, where the  measure $\nu$ is absolutely continuous with respect to Lebesgue measure, with a density $\cL (\di \mu)$.

Stieltjes transforms $f$ have the property that $-f$ is reciprocally convex (in terminology introduced in \cite{recko}, a function $g(x)$ is
reciprocally convex if  it is defined for $x>0$ and concave there, whereas $g(1/x)$ is convex).  As proved in \cite{recko},
each reciprocally convex function  generates
an increasing  sequence of quasi-arithmetic means, and hence  \cm functions that are also Stieltjes
transforms are interesting as a tool for generating  means.

}

\eitem

\qitem{Completely monotone probability densities}{ Let $f$ be a probability density with respect to Lebesgue measure on $[0,+\infty)$, that is,
\[ \int_0^{+\infty} f(x)\di x  = 1 \qquad \mbox{and \quad $f(x)\geq 0$ for all $x\geq 0$}. \]

Then $f$ is a \cmf if and only if
(\ref{laptra}) holds, which, after integration with respect to $x\in (0,+\infty)$ gives (via Fubini theorem for $f\geq 0$)
\[ 1 = \int_0^{+\infty} \frac{1}{t} \di \mu (t) .\]
Defining a new probability measure $\nu$ by $\nu (B)=\int_B \frac{1}{t}\di \mu (t)$, we have that
\be
\label{pdcm1}
f(x)= \int_{[0,+\infty)} t e^{-xt}\di \nu (t) = \int_{[0,+\infty)} t e^{-xt}\di G (t),
\ee
where  $G$ is the distribution function for $\nu$.
The function $x\mapsto te^{-xt}$ is the density of exponential distribution $\expd (t)$. Therefore,
a density $f$ of a probability measure on $(0,+\infty)$ is a \cmf if and only if it is a mixture of exponential
densities. Note that (\ref{pdcm1}) can be written as $f(x)= \E (Te^{xT})$, where $T$ is a random variable with distribution function $G$; by
letting $S=1/T$ we find that
\be
\label{pdcm2}
f(x)=\E \left( \frac{1}{S} e^{\frac{x}{S}}\right) = \int_{[0,+\infty)} \frac{1}{s} e^{-x/s} \di H(s),
\ee
where $H$ is the distribution function of $S$. The latter form is taken as a definition of what is
meant by a \cm density in \cite[18.B.5]{marolk}; this is more natural than (\ref{pdcm1}) because
the mixing measure $H$ is defined on values of expectations ($s$) of exponential distributions in the mixture, rather then
on their reciprocal values as in (\ref{pdcm1}).
}

\eitem

\section{Further properties and connection with infinitely divisible measures}
\label{proprel}

Starting from the mentioned representations of \cm functions,
an interesting criterion for equality of two \cm functions is derived in \cite{fell68}:

\qitem{Lemma}{\em \label{eqcm} If
$f$ and $g$ are \cm functions and if $f(x_n)=g(x_n)$  for a  positive sequence $\{x_n\}$ such that the series  $\sum_n 1/x_n$ diverges, then
$f(x)=g(x)$ for all $x\geq 0$.}

As a corollary to Lemma \ref{eqcm}, we can see that if \cm functions
$f$ and $g$ agree in any subinterval of $(0,+\infty)$, then $f(x)=g(x)$ for all $x\geq 0$.
A converse result, which  is also proved in \cite{fell68} is more surprising: If $f$ is \cm and if the series $\sum_n 1/x_n$ converges,
then there exists another \cm
function $g\neq f$, such that $f(x_n)\ne g(x_n)$ for all $n$.

\eitem

\qitem{Convolution and infinitely divisible measures}{\label{cidm}
 Given measures $\mu$ and $\nu$ on $[0,+\infty)$ and their distribution functions $g_{\mu}$ and $g_{\nu}$, we define
the convolution $\mu \ast \nu$ as a measure with the distribution function defined by
\be
\label{convm}
g_{\mu \ast \nu} (t) = \int_{[0,t]} g_{\mu}(t-u) \di g_{\nu} (u) = \int_{[0,t]}  g_{\nu}(t-v) \di g_{\mu} (v)
\ee
To show equality of integrals above, we use the formula for integration by parts in Lebesgue-Stieltjes integral
(see \cite{hew60} or \cite{lesti}) and
note that the function $u\mapsto g_{\mu}(t-u)$ is continuous from the left, while $u\mapsto g_{\nu}(u)$ is continuous
from the right, hence the additional term due to discontinuities in the integration by parts formula equals zero, that is,
\[ \int_{[0,t]} g_{\mu}(t-u) \di g_{\nu} (u) = -\int_{[0,t]}  g_{\nu}(u) \di g_{\mu} (t-u) \]
and then we apply change of variables in the last integral, $u=t-v$.

Repeated convolution is defined by induction, using associativity. In particular, the $n$th convolution power of a measure $\mu$, denoted by
$\mu^{n\ast}$ is defined by $n-1$ repeated convolutions $\mu\ast\mu \ast \cdots \ast\mu$.

A measure $\mu$ is called infinitely divisible (\id) if for every natural number $n$ there exists a measure $\mu_n$ such that
$\mu=\mu_n^{n\ast}$.

}

\eitem

In the next two lemmas we collect some basic properties of \cm functions. For a collection of other properties we refer to \cite{milsam01}.

\medskip

\qitem{Lemma}{\em \label{basprop1} If $f$ and $g$ are \cm functions with $f=\cL (\di \mu)$ and $g=\cL (\di \nu)$, then for $a>0$,
\[ af= \cL (d(a\mu)),\quad f+g  =\cL (\di (\mu +\nu)),\quad fg = \cL (\di (\mu\ast \nu)) .\]
Therefore, if $f,g$ are \cm then $af+bg$ ($a,b>0$) and $fg$ are also \cmt }

{\bf Proof. } First two properties follow from the definition of Laplace transform. The third property for arbitrary positive measures
is proved in \cite[p. 434]{fell2}.
\eitem

\qitem{Lemma}{\em \label{basprop2} (i) If $g'$ is  \cmz then the function  $x\mapsto f(x)=e^{-g(x)}$ is  \cmt

(ii) If $\log  f$ is  \cmz then $f$ is \cm (the converse is not true).

(iii) If $f$ is \cm and $g$ is a positive function with a \cm derivative,
then $x\mapsto f(g(x))$ is \cmt
}

{\bf Proof. } To prove (i), let $h(x)=e^{-g(x)}$ and note that $h>0$ and $h'=-g'h<0$. Then by induction, using Leibniz chain rule, it follows
that $(-1)^n h^{(n)} >0$. In particular, if $\log f$ is  \cmz  then $(-\log f)'$ is also \cmz and (ii) follows from (i) with $g=-\log f$.
The function $x\mapsto e^{-x}$ is a \cm function but its logarithm is not the one, so the converse does not hold.
For (iii), we note that $f=\cL (\di\mu)$ for some positive measure $\mu$, hence
\be
\label{derfg}
\der{}{x} f(g(x)) = -g'(x)\int_0^{+\infty} e^{-g(x)t} t\di \mu(t)
\ee
By part (i), the function $x\mapsto e^{-g(x)t}$ is \cm  for every $t>0$, and so the function $x\mapsto g'(x)e^{-g(x)t}$ is
also \cm as a product of two \cm functions. Then from representation (\ref{derfg}) it follows that the first derivative of $-f(g(x))$ is \cmz
which together with positivity of $f$ yields the desired assertion. \eop

Note that if we can find measures $\nu$ and $\nu_t$ in representations $g'(x)=\cL (\di \nu)$ and $e^{-g(x)t}=\cL(\di \nu_t)$, then
from (\ref{derfg}) we find that
\be
\label{derfga}
\der{}{x} f(g(x)) = -\int_0^{+\infty} t \int_0^{+\infty} e^{-ux}\di (\nu \ast \nu_t) (u) \di \mu(t).
\ee

\eitem

It turns out that \cm functions $f$ of the form as in (i) of Lemma \ref{basprop2} are Laplace transforms of \id\
 measures. If $f(0)=1$, the associated measure is a probability measure, which is the case that is of interest in
applications. Proofs of statements of the next lemma can be found in \cite{fell2}.

\medskip

\qitem{Lemma}{\em \label{ltid} (i) A function $f$ is the Laplace transform of an \id\  probability measure if and only if
\be
\label{ltid1}
f(x)=e^{-g(x)},
\ee
where $g$ is a positive function with a \cm derivative and $g(0)=0$.  Equivalently, $f$ is the Laplace transform of an id positive measure if
and only if $f(x)>0$ for all $x>0$,  and the function $x\mapsto -\log f(x)$ has a \cm derivative. This measure is a probability measure
if and only if $f(0_{+})=1$.

(ii) A function $f$ is the the Laplace transform of an \id\ probability measure if and only if
\be
\label{ltid2}
-\log f(x) = \int_0^{+\infty} \frac{1-e^{-xt}}{t}\di \mu (t),
\ee
where $\mu$ is a positive measure such that
\be
\label{ltid3}
\int_1^{+\infty} \frac{1}{t}\di\mu(t) <+\infty.
\ee
}

\medskip

\eitem

\qitem{Remarks}{$1^{\circ}$ If $\log f$ is  \cmz then $-\log f$ has a \cm derivative and by Lemma \ref{ltid}(i),
 $f = \cL (\di \mu)$, where $\mu$ is an \id\ positive measure. By \cm property of $\log f$, we have that $\log f = \cL (\di \nu)$,
 where $\nu$ is some other positive measure. Note that positivity of $\nu$ implies that $\log f(0)>0$, that is,
 $\mu ([0,+\infty))=f(0)>1$ and so, $\mu$ can not be a probability measure.

$2^{\circ}$. Non-negative  functions with a \cm first derivative have a special name - {\em Bernstein functions};
Lemmas \ref{basprop2} and \ref{ltid} explain their role in probability theory; more about this class of functions can
be found in \cite{bernf}.

}

\eitem

\section{Majorization, convexity and logarithmic convexity}
\label{majocolo}

A good source for studying all three  topics that are very much interlaced, is the book \cite{marolk}. In this short digest we
include only necessary definitions and results that one can need for understanding a connection with \cm functions.

\qitem{Majorization and Schur-convexity}{\label{major}   For a vector  $\bx\in \R^n$ define $x_{[i]}$ to be the $i$th largest coordinate of $\bx$, so that
\[ x_{[1]} \geq x_{[2]}\geq\cdots\geq  x_{[n]}.\]
We say that $\bx$ is majorized by $\by$ in notation $\bx \prec \by$ if
\[ \sum_{i=1}^k x_{[i]} \leq \sum_{i=1}^k y_{[i]} \quad {\rm for }\ k=1,2,\ldots, n-1\quad {\rm and}\quad
\sum_{i=1}^n x_{[i]} = \sum_{i=1}^n y_{[i]} \]

For example, $(1,1,1)\prec (2,1,0)$. Clearly, majorization is invariant to permutations of coordinates of vectors.

A function $f$ which is defined on a symmetric set $S\subset \R^n$ ($S$ is symmetric if $\bx\in S$ implies that $\by \in S$ where
$\by$ is any vector obtained by permuting the coordinates of $\bx$) is called Schur-convex if for any $\bx,\by \in S$,
\be
\label{scdef}
\bx\prec \by \implies f(\bx) \leq f(\by).
\ee

The following result, due to A. M. Fink\cite{fink} reveals an interesting relationship between concepts of  Schur-convexity
and complete monotonicity.

}
\medskip

\qitem{Lemma}{\em \label{lemfink} For a \cm function $f$ and a non-negative integer vector $\bm = (m_1,m_2,\ldots, m_d)$ of a dimension $d>1$.
 let
\[ u_x (\bm) =(-1)^{m_1}f^{(m_1)}(x)(-1)^{m_2}f^{(m_2)}(x)\cdots (-1)^{m_d}f^{(m_d)}(x) .\]
Then $u_x (\bm)$ is a Schur-convex function on $\bm$ for every $x>0$ and $d>1$.

}

An important corollary of \ref{lemfink} is  with $d=2$, taking $\bm = (1,1)$ and $\bn =(2,0)$. Clearly, $\bm \prec\bn$ and from the above
definition of Schur-convexity we get that $u_x(1,1) \leq u_x(0,2)$, that is, $(f'(x))^2 \leq f(x)f''(x)$, which is, knowing that
$f(x)>0$, equivalent to $(\log f(x))''\geq 0$. We formulate this result as a separate lemma.

\eitem

\qitem{Lemma}{\em \label{cflogc} Any \cm function $f$ is log-convex, i.e., the function $\log f(x)$ is convex.}

\eitem

A converse does not hold, for example  the Gamma function restricted to $(0,+\infty)$ is log-convex, but it is not \cm. However,
the fact that each \cm function is also log-convex, helps us to search for possible candidates for complete monotonicity only
among functions that are log-convex. In addition, there is a very rich theory that produces inequalities using convexity or Schur-convexity,
and we can use it for \cm functions.

Log-convexity of \cm  functions is equivalent to decreasing of the ratio $f'(x)/f(x)$, and (arguing that $f^{(2k)}$ and $-f^{(2k+1)}$ are
\cmn) this implies

\qitem{Corollary}{\em \label{decrat} If $f$ is a \cm function, then the ratio
\[ x\mapsto \left|\frac{f^{(k+j)}(x)}{f^{(k)}(x)}\right|  \]
is decreasing for every integers $k,j$. }

\eitem
In the next lemma we give two consequences of convexity and log-convexity of \cm functions. Similar inequalities for \cm functions can be found in
\cite{kimber74}, but with more involved proofs.

\qitem{Lemma}{\label{940}\em If $f$ is completely monotone, then

\beq
\label{supad}
 f(x)+f(y)  \leq & f(x-\eps) + f(y+\eps)&\leq f(0)+f(x+y), \label{supad1}\\
 f(x)f(y)   \leq & f(x-\eps)f(y+\eps) &\leq  f(0)f(x+y) \label{supad2}
 \eeq
 where $0\leq \eps<x<y$, assuming that $f(0)$ is defined as $f(x_{+})$ (as in \ref{limpro}, finite or not).

}

{\bf Proof. } If $\fai$ is a convex function, then the divided difference
\[ \Delta_{\fai,\eps} (x) = \frac{ \fai (x) - \fai (x-\eps)}{\eps} \]
is increasing with $x$, hence in the present setup, $\Delta_{f,\eps} (x)\leq \Delta_{f,\eps} (y+\eps)$ and $\Delta_{f,x-\eps} (x-\eps)
\leq \Delta_{f,x-\eps} (x+y) $,
which proves (\ref{supad1}). The same proof holds for (\ref{supad2}), but with $\log f$ in place of $f$. \eop

Let us note that under assumptions of Lemma \ref{940}, $(x,y)\prec (x-\eps,y+\eps)\prec (0,x+y)$, and so we have just proved
that the functions $(x,y)\mapsto f(x)+f(y)$ and $(x,y)\mapsto f(x)f(y)$ are Schur-convex on $\R_{+}\times \R_{+}$. More generally,
for any $f$ being \cm, the functions of $n$ variables
\be
\label{sconv}
\sum_{i=1}^n f(x_i)\quad \mbox{and}\quad \prod_{i=1}^n f(x_i)
\ee
are Schur-convex on $\R_{+}^n$.  For a proof of this statement see \cite{marolk}.

\eitem

Finally,  the fact that $f'$ is concave (i.e, $f'''<0$) is equivalent to  each of  three inequalities in the next lemma\cite{code,repje}.

\qitem{Lemma}{\label{146} \it For a \cm function $f$, it holds
\be
\label{conc1}
\frac{f'(x)+f'(y)}{2}< \frac{f(y)-f(x)}{y-x}< f'\left(\frac{x+y}{2}\right),\quad \mbox{for all $x,y >0$}
\ee

\be
\label{conc2}
\frac{f(y)-f(x)}{y-x} < \frac{f(y-\eps)-f(x+\eps)}{y-x-2\eps},\quad \mbox{for  $0<x<y$ and $0<\eps<\frac{y-x}{2}$}
\ee

}

\section{Inversion formulae}
\label{inver}

It is sometimes easier to find a measure $\mu$ that corresponds to function $f$ via Laplace transform
in (\ref{laptras}) then to show that $f$ is  \cm by verifying the definition; in view of applications, it is definitely useful and desirable
to know the associated measure.   In many cases we can use properties of Laplace transform and
the tables that can be found in textbooks. In many applications the Laplace transform is not limited
to real argument, and it is more common to define $f(z)$ by (\ref{laptras}),
where complex argument $z$ belongs to some half space $\Re z \geq a$, for some positive $a$. We may use the power of  complex Laplace
transform calculus applied to real function of real argument, due to well known properties of regular functions.

Due to similarity between Fourier transform, complex Laplace transform and real Laplace transform,  we
may use inversion formulae for all three mentioned classes, whenever it is appropriate.
 In probability theory, for a random variable $Z$, the function
$x\mapsto \E e^{ixZ}$ (which corresponds to Fourier transform, except the sign in the exponent) is called the
characteristic function, whereas the  real Laplace transform (mind the sign!) $x\mapsto \E e^{xZ}$ is called the moment generating function.
There are several formulas that can be found in textbooks,  but  we
will mention here only a not widely known inversion theorem that enables finding a finite measure
$\mu$ defined on Borel sets of $\R$, provided that we know its characteristic function
\be
\label{chfm}
\fai (x) = \int_{-\infty}^{+\infty} e^{itx} \di F(t),
\ee
where $F(t) = \mu \{ (-\infty,t]\}$. The following result (given here in a slightly generalized version) is due to Gil-Pelaez \cite{gilp51}.

\qitem{Lemma}{\em \label{fuinv} For $\fai$ and $F$ as in (\ref{chfm}), with $\fai (0)$ being finite,  we have that, for all $t\in \R$,
\be
\label{fuinv1}
\frac{ F(t)+F(t_{-})}{2} = \frac{\fai(0)}{2} - \frac{1}{\pi}\int_0^{+\infty}\Re \left( \frac{e^{-itx}\fai(x)}{ix}\right) \di x.
\ee
 }
Note that the underlying measure here need not necessarily be restricted to the positive part of the real axis.
As an example of how (\ref{fuinv1})
can be used to determine a measure $\mu$ such that  $f=\cL (\mu)$,
consider a simple case $f(x)=e^{-x}$, where we already know that the measure is Dirac at $t=1$.
 Supposing  that we wish to use (\ref{fuinv1}) to derive this, note that if $f$ is the Laplace transform of $\mu$,
then its characteristic function is $\fai (x)=
f(-ix)=e^{ix}$, and  (\ref{fuinv1}) yields (assuming that $t$ is a point of continuity of $F$)
\be
\label{fuinv2}
F(t)= \frac{1}{2}- \frac{1}{\pi} \int_0^{+\infty} \frac{\sin x(1-t)}{x}\di x .
 \ee
Knowing that
\[ \int_0^{+\infty}\frac{\sin (ax)}{x} \di x = \frac{\pi}{2}\sign a , \]
we find that $F(t)= 0$ for $t<1$ and $F(t)=1$ for $t>1$, hence (by right-continuity and non-decreasing of $F$),
the corresponding measure $\mu$ is indeed a Dirac measure at $t=1$.

\eitem

For other formulas and methods, including numerical evaluation of inverse, see \cite{cohen07}. In the next lemma we
complement some examples from \cite{milsam01} by effectively finding
the corresponding measure.

\qitem{Lemma}{\label{milsamplus} We have the following representations:
\be
\label{milsam1}
e^{-ax} = \cL (\di \delta_a (t)),
\ee
where $\delta_a$ is the  probability measure with unit mass (Dirac measure)  at $a\geq 0$\ ;

\be
\label{milsam2}
\frac{1}{(ax+b)^c} = \cL \left(e^{-bt/a}\frac{t^{c-1}}{a^c\Gamma(c)}\right),\quad a,b,c\geq 0, a^2+b^2>0\ ;
\ee

\be
\label{milsam3}
\log\left(a+\frac{b}{x}\right) = \cL (\di \mu(t)) \qquad a\geq 1, b>0,
\ee
where the measure $\mu$ is determined by its distribution function
\[ \mu ([0,t]) = \log a + \int_0^x \frac{1-e^{-bs/a}}{s}\di s ;\]

\be
\label{milsam4}
\frac{\log (1+x)}{x} = \cL ( E_1(t)) ,
\ee
where (see \cite[p.56]{absteg}) $E_1$ is exponential integral
\[E_1(t)= \int_1^{+\infty}e^{-tu}\frac{\di u}{u} ;\]

\be
\label{milsam5}
e^{a/x} = \cL \left(\di \delta_0 (t) + \frac{aI_1(2\sqrt{at})}{\sqrt{at}}\di t\right),
\ee
where $I_1$ is a modified Bessel function as defined in \cite{absteg}.
}

{\bf Proof. } The relation (\ref{milsam1}) is obvious, and (\ref{milsam2}) is a consequence of standard rules for (complex) Laplace transform:
\[
\cL \left(e^{-bt/a}\frac{t^{c-1}}{a^c\Gamma(c)}\right) = \frac{1}{a^c\Gamma (c)}\cL (t^{c-1})(x-b/a)
= \frac{1}{a^c\Gamma (c)} \cdot \frac{\Gamma (c)}{(x-b/a)^c}=\frac{1}{(ax+b)^c}
 \]
 To prove (\ref{milsam3}), denote its left side  by $f$, and observe that, by (\ref{milsam2}),
 \[ f'(x) = \frac{a}{ax+b}-\frac{1}{x} = \cL \left(e^{-bt/a}-1\right).\]
 Now we use the rule
 \[ \cL \left(\frac{g(t)}{t}\right) = \frac{1}{x}\int_x^{+\infty} \cL (g(t))[y]\di y \]
 to conclude that
 \[ f(x) = \log a - \int_x^{+\infty} f'(y)\di y = \cL (\log a \di \delta_0) - \cL \left(\frac{e^{-bt/a}-1}{t}\right),\]
 which yields (\ref{milsam3}).
 To prove (\ref{milsam5}), we note that
 \[ \frac{a^k}{k!x^k}= \cL\left(\frac{a^k t^{k-1}}{k!(k-1)!}\right),\]
which  tells us that
 \[ e^{\frac{a}{x}} = 1+\cL\left( \sum_{k=1}^{+\infty} \frac{a^k t^{k-1}}{(k-1)!k!}\right). \]
 Now we observe that
 \[\sum_{k=1}^{+\infty} \frac{a^k t^{k-1}}{(k-1)!k!} =\frac{aI_1(2\sqrt{at})}{\sqrt{at}}, \]
 and (\ref{milsam5}) follows.

The simplest way to prove (\ref{milsam4}) would be to perform an integration on the right hand side and show that it yields the left side.
However, in order to show the derivation, we  start with the observation that
\[ f(x):=\frac{\log(1+x)}{x}= F(1,1,2;-x), \]
where $F(a,b,c; \cdot)=\ _2F_1(a,b,c;\cdot)$ is Gauss' hypergeometric function;
hence there is the following integral representation \cite{absteg}:
\[ f(x) = \int_0^1 \frac{\di s}{1+sx} .\]
Now we use (\ref{milsam3}) to find that
\[ \frac{1}{1+sx}=\frac{1}{s} \int_0^{+\infty} e^{-xt}e^{-t/s}\di t, \]
and, exchanging the order of integration, we find that
\[ f(x) = \int_0^{+\infty} e^{-xt}\left(\int_0^1 e^{-t/s}\frac{\di s}{s}\right)\di t. \]
Finally, a change of variables $1/s=u$ in the inner integral shows that it is equal to $E_1(t)$, and the formula is proved.

\eitem

\section{Some examples related to the Gamma function}
\label{examp}

\medskip

Functions related to the Gamma function are good candidates to be \cm, and there is a plenty of such results in literature. The  function
$g(x)=\log \Gamma (x)$ is a unique convex solution of the Krull's functional equation
\be
\label{krulle}
g(x+1)-g(x) = f(x),\qquad x>0,
\ee
with $f(x)=\log x$ and with $g(1)=0$. The same equation, but with $f(x)=(\log x)^{(n+1)}$, $n=0,1,2,\ldots$ has for its solutions
functions $\Psi^{(n)}(x)=(\log \Gamma (x))^{(n+1)}$.
Although $\log x $ is not \cmn, all its derivatives are monotone functions, which automatically
implies the same property for $\Psi^{(n)}(x), n\geq 2$ and alike functions via the following result (see \cite{cnvx07}).

\qitem{Lemma}{\label{krullcm} \it Suppose that $x\mapsto f(x)$
is a function of the class $C^{\infty}(0,+\infty)$ with all derivatives being monotone functions, with $f'(x)\tends 0$ as $x\tends +\infty$.
 Then there is
 a unique (up to an additive constant)  solution $g$ of (\ref{krulle}) in the class $C^{\infty}$, with

\be
\label{fdercm}
 g'(x)      = \lim_{n\tends +\infty} \left( f(x+n) - \sum_{k=0}^n f'(x+k) \right)
 \ee

 and

 \be
 \label{odercm}
g^{(j)}(x) = -\sum_{k=0}^{+\infty} f^{(j)} (x+k)\qquad (j\geq 2).
\ee

}
From (\ref{fdercm}) and (\ref{odercm}) it follows that, if $\pm f$ is  \cm (or if only $\pm f''$ is such),
then $\mp g''$ is a \cmn, while  $\pm g$ and $\pm g'$ need not be \cmn .
Our first example is formulated in the form of a lemma, and its proof provides a pattern that can be used in many similar cases.

\qitem{Lemma}{\label{lem0} \it The function
\[ W(x)=-(\log \Gamma (x)-(x-1)\log x)''=\frac{1}{x}+\frac{1}{x^2}-\Psi'(x)  \]
has the following integral representation
\be
\label{wi}
 W(x)=\int_0^{+\infty} \left(1+t-\frac{t}{1-e^{-t}}\right)e^{-xt}\di t
 \ee
and it is a \cmf.}

{\bf Proof. } The integral representation follows from
\be
\label{psin}
\cL \left( (-1)^{n+1}\frac{t^n}{1-e^{-t}}\right)=\Psi^{(n)}(x), \qquad n=1,2,\ldots,
\ee
and
\be
\label{recxn}
\cL (t^a) = \frac{\Gamma (a+1)}{x^{a+1}},\qquad a>-1.
\ee
The \cm property follows from positivity of the function under integral sign, which is equivalent to the inequality $e^t>1+t$ for $t>0$.

\eitem

{\bf Remark. } The function $g(x)=\log \Gamma (x) - (x-1)\log x$ satisfies the functional equation (\ref{krulle}) with
$f(x)=\log\left(\frac{x}{x+1}\right)^{x}$; it can be easily checked that $f''$ is \cmn, hence from Lemma \ref{krullcm} we can conclude without
any additional work that $-W$ as defined above is \cm.

\qitem{Example}{\label{exa200} \it For $a\geq 0$ and $x>0$, let
\[ G_a(x)= \log \Gamma (x) -\left( x-\frac{1}{2}\right) \log x -\frac{1}{12}\Psi'(x+a) +x-\frac{1}{2}\log (2\pi).\]
The following  representation holds:
\be
\label{gaga}
G_a(x)=\int_0^{+\infty} \frac{t-2+(2+t)e^{-t} -(t^3/6) e^{-at}}{2t^2(1-e^{-t})}e^{-xt} \di t.
\ee
The function $x\mapsto G_a(x)$  is \cm if and only if $a\geq 1/2$ and the function $x\mapsto - G_a(x)$ is \cm if and only if $a=0$.}

{\bf Proof. } Starting with
\[  G_a''(x)= \Psi'(x)-\frac{1}{12}\Psi'''(x+a)-\frac{2}{x}+\frac{x-1/2}{x^2}, \]
it easy to show (in a similar way as in Lemma \ref{lem0}) that
\be
\label{fair}
G_a''(x)=\int_0^{+\infty} \frac{t-2+(2+t)e^{-t} -(t^3/6) e^{-at}}{2(1-e^{-t})}e^{-xt} \di t.
\ee
Further, we have that
\[ \lim_{x\tends +\infty} G_a(x)=\lim_{x\tends +\infty} G_a'(x) = 0,\]
and
\[ G_a(x)=\int_x^{+\infty}\int_v^{+\infty}G_a''(u)\di u\di v ,\]
hence (\ref{gaga}) holds.
 The complete
monotonicity is related to the sign of the function
\be
\label{ha}
 h_a (t) = t-2+(2+t)e^{-t} -\frac{t^3}{6} e^{-at} .
 \ee
The function $G_a$ is \cm if and only if $h_a(t)\geq 0$ for all $t\geq 0$. From (\ref{ha}) we see that this is equivalent to
\be
\label{haeq}
 a\geq \frac{\log 6+\log ((2+t)e^{-t} + t -2)-3\log t }{-t}:=u(t)
 \ee
Using standard methods, we can find that $u$ is a decreasing function, hence
\[ u(t) \leq \lim_{t\tends 0_{+}} u(t) = \frac{1}{2}, \]
and so, (\ref{haeq}) holds if and only if $a\geq 1/2$.

Further, $-G_a$ is \cm if and only if $h_a(t) \leq 0$ for all $t\geq 0$, which is equivalent to
\be
\label{haeqa}
a\leq u(t),
\ee
where $u(t)$ is defined in (\ref{haeq}). Since $u$ is decreasing, we have that
\[ u(t) \geq \lim_{t\tends +\infty} u(t) =0 ,\]
and so, (\ref{haeqa}) holds if and only if $a\leq 0$, that is, $a=0$.

{\bf Remark.}  Let
\[ F_a(x) = \log \Gamma (x) -\left( x-\frac{1}{2}\right) \log x -\frac{1}{12}\Psi'(x+a), \qquad a\geq 0, x>0\]

This function is studied in \cite[Theorem 1]{digam}, where it is shown that  $x\mapsto F_0(x)$ is  concave on $x>0$ and that
$x\mapsto F_a(x)$ is convex on $x>0$ for $a\geq \frac{1}{2}$. Since $F_a(x)''=G_a''(x)$ where $G_a$ is defined as above, this example gives much stronger
statement.

\eitem

\qitem{Example}{\label{exa300} \it For $b\geq 0$ and $c\geq 0$, let
\be
\label{fun301}
f_{b,c}(x) = \frac{e^x \Gamma (x+b)}{x^{x+c}},\qquad x>0.
\ee
The function
\be
\label{fun3015}
\fai_{b,c}(x)= \log f_{b,c}(x) = x+\log \Gamma (x+b) - (x+c)\log x
\ee
 is \cm if and only if   $b\geq \frac{1}{2}+\frac{1}{\sqrt{12}}$ and $c=b-\frac{1}{2}$ and then it has the representation
\be
\label{fun302}
\fai_{b,b-\frac{1}{2}}(x)
= \int_{[0,+\infty)}\frac{1}{t^2}\left( \frac{t e^{-bt}}{1-e^{-t}}+t\left(b-\frac{1}{2}\right)-1\right)\di t. \qquad x>0.
\ee

}

{\bf Proof. } By expanding  $\log\Gamma (x+b)$ in (\ref{fun3015}) by means of  Stirling's formula \cite[p.258]{absteg}, it
follows that, for $\delta=b-c\neq \frac{1}{2}$,
\[ \lim_{x\tends +\infty} \fai_{b,c} (x) = \left(\delta -\frac{1}{2}\right) \cdot (+\infty),  \]
so $\fai_{b,c}$ is not a \cm function (see \ref{limpro}). Let $\delta = 1/2$ and let
\[ G_b (x):= \fai_{b,b-\frac{1}{2}}(x)= x+\log \Gamma (x+b) - \left(x+b-\frac{1}{2}\right)\log x.\]
Further, we find without difficulties that
\be
\label{fun3025}
 \lim_{x\tends +\infty} G_b (x) = \lim_{x\tends +\infty} G_b'(x) =0
 \ee
and that
\[ G_b'' (x) = \Psi'(x+b)-\frac{1}{x} +\frac{b-\frac{1}{2}}{x^2} .\]
In the same way as shown in Lemma \ref{lem0}, we find that  $G_b''(x)= \cL (h_b(t)\di t)$, where
\be
\label{fun3027}
h_b(t) = \frac{t e^{-bt}}{1-e^{-t}}+t\left(b-\frac{1}{2}\right)-1 .
\ee

By standard methods we find that
\be
\label{fun303}
 h(t)=\left( \frac{b^2}{2} -\frac{b}{2} +\frac{1}{12}\right) t^2 + o(t^2)\qquad (t\tends 0),
 \ee
so the Laplace transform $G_b(x)$ of the function $t\mapsto g(t)/t^2$ exists for all $x>0$ and applying Fubini theorem as
in Example \ref{exa200} and using (\ref{fun3025}) we find that
\[ G_b (x) = \int_x^{+\infty}\int_v^{+\infty}G_a''(u)\di u\di v = \int_{[0,+\infty)} \frac{h(t)}{t^2}e^{-tx}\di t ,\]
which is the representation (\ref{fun302}). Then $G_b$ will be \cm if and only if if $h(t)\geq 0$ for each $t\geq 0$ (see
Remark \ref{rmden}). By (\ref{fun303}) we have that $h(0)<0$ for $b\in (b_1,b_2)$, where $b_{1,2}= \frac{1}{2}\pm \frac{1}{\sqrt{12}}$;
further, $c=b-1/2>0$ gives $b>1/2$, so only $b\geq b_2$ remains as a possibility. It is straightforward to check that $\pder{h_b(t)}{b}>0$
for all $t\geq 0$, so it suffices to show that $h_{b_2}(t)\geq 0$ for $t\geq 0$, which can be done along the lines of \cite{guosri08}.

{\bf Remark. } Complete monotonicity of functions $f_{b,c}$ and $\fai_{b,c}$ for various values of parameters
was discussed in \cite{guosri08} and \cite{guosri07}.
Let us remark that, by Lemma \ref{basprop2}, the function $f_{b,c}$ is \cm whenever $\fai_{b,c}$ is the one. \eop

\medskip

Let us mention that the Barnes function $G(x)$ satisfies the relation
\[ \log G(x+1) - \log G(x) = \log \Gamma (x),\quad x>0, \]
which is (\ref{krulle}) with $g=\log \Gamma$. Here also the function
$x\mapsto (\log G(x))''= 2\Psi'(x) +(x-1)\Psi''(x)$ is \cmn.
More details about the properties of the $G$-function as a solution of Krull's equation can be found in \cite{mimo11}.
\eitem

\medskip

{\small \sc
\noindent University of Belgrade\\
Faculty of Electrical Engineering\\
P.O. Box 35-54, 11000 Belgrade, Serbia}\\
emerkle@etf.rs

\smallskip

\noindent{\small\sc and

\smallskip

\noindent Union University\\
Ra\v cunarski fakultet\\
Kneza Mihaila 6, 11000 Belgrade, Serbia

}

\end{document}